\documentclass[final,3p,times]{elsarticle}

\usepackage{amssymb}

\usepackage{graphicx}
\usepackage{amssymb}
\usepackage{amsmath}
\usepackage{cases}
\usepackage{amsfonts}
\usepackage{epsfig}
\usepackage{booktabs}
\usepackage{upgreek}
\usepackage{graphicx,tabularx,multirow, color}
\usepackage[active]{srcltx}
\makeatletter

\newcommand{\R}{\mathbb{R}}

\newcommand{\diff}{\,\mathrm{d}}
\newcommand{\bbx}{\mathbf{x}}
\newcommand{\bby}{\mathbf{y}}

\newcommand{\bbA}{\mathbf{A}}

\newcommand{\bbr}{\mathbf{r}}

\newcommand{\bbs}{\mathbf{s}}
\newcommand{\bbb}{\mathbf{b}}
\newcommand{\bba}{\mathbf{a}}

\newcommand{\bbf}{\mathbf{f}}

\newcommand{\im}{\mathrm{i}}

\makeatother

\begin{document}

\begin{frontmatter}



\title{An efficient adaptive frequency sampling scheme for large-scale transient boundary element analysis}


\author[a1,a11]{Jinyou Xiao\corref{cor1}}
\ead{xiaojy@nwpu.edu.cn}
\cortext[cor1]{Corresponding author}
\address[a1]{School of Astronautics, Northwestern Polytechnical University, Xi'an 710072, China}
\address[a11]{Institute for Computational Mechanics and Its Applications, Northwestern Polytechnical University, Xi'an 710072, China}

\author[a2]{Junjie Rong}
\author[a2]{Wenjing Ye}
\ead{mewye@ust.hk}
\address[a2]{Department of Mechanical and Aerospace Engineering, Hong Kong University of Science and Technology, Hong Kong}

\author[a3]{Chuanzeng Zhang}
\ead{c.zhang@uni-siegen.de}
\address[a3]{Department of Civil Engineering, University of Siegen, D-57068 Siegen, Germany}

\begin{abstract}
The frequency-domain approach (FDA) to transient analysis of the boundary element method, although is appealing for engineering applications, is computationally expensive. This paper proposes a novel adaptive frequency sampling (AFS) algorithm to reduce the computational time of the FDA by effectively reducing the number $N_{\rm{c}}$ of sampling frequencies.
The AFS starts with a few initial frequencies and automatically determines the subsequent sampling frequencies. It can reduce $N_{\rm{c}}$ by more than 2 times while still preserving good accuracy. In a porous solid model with around 0.3 million unknowns, 4 times reduction of $N_{\rm{c}}$ and the total computational time is successfully achieved.

\end{abstract}


\end{frontmatter}


\section{Introduction}

The boundary element method (BEM) has been extensively used to solve elastodynamic problems in many fields of science and engineering, including fracture mechanics \cite{SS87,Zh93}, seismology \cite{KB84,CBS09} and soil-structure interactions \cite{PG99}. This paper is devoted to the efficient solution of transient elastodynamic problems.

According to the different solution strategies in time space, the BEM for treating such problems generally follows two approaches, namely, time-domain approaches and frequency-domain approaches; see, e.g., the reviews by Beskos \cite{Bde97} and
Costabel \cite{CM04}.
Time-domain approaches can be further classified into time-stepping methods and the space-time integral equation method. In these methods the physical problems are directly solved in the real time domain, thus one can observe the phenomenon as it evolves.
However, such methods require an adequate choice of the time step size. An improper
time step could lead to instability or numerical damping. For recent development of time-domain methods, see e.g. \cite{SA97,ADF12}.
Beside these methods there exist the possibility to solve the time-domain boundary integral equation with the so-called convolution quadrature method proposed by Lubich \cite{Lubich88}, which provides a straightforward way to obtain a stable time-stepping scheme using the Laplace transform of the kernel
function \cite{SA97,Zh2000,BMS12}.

The frequency-domain approach based on the Laplace transforms offers another attractive approach for transient analysis \cite{AM87,PTB05,PGS10}. In this approach, one solves the frequency-domain boundary integral equations at a series of discrete frequencies, then obtains the time-domain responses by employing certain numerical inverse Laplace transform  methods \cite{K13}. The Fourier series method (FSM) \cite{Crump76,Cohen07} is one of the most popular methods in computing the inverse Laplace transform which has found wide applications in boundary element transient analysis, see e.g. \cite{MB81,GaulS99,XYCZ12,XYW13}. It works by truncating the infinite Bromwich contour integral of inverse Laplace transform into a finite one, and evaluating the finite integral using the trapezoidal rule based on equally-spaced integration points (i.e., sampling frequencies). As such, the time-domain responses can be efficiently obtained by using the fast Fourier transform \cite{MR08}; see Section \ref{S-prob-fda}. In real applications the number of sampling frequencies in the FSM can often be more than one hundred. In large-scale BEM analysis this undoubtedly implies, a quite huge, if not prohibitive, computational burden.

This paper is devoted to the effective reduction of the number of sampling frequencies in the frequency-domain approach. The outcome is an adaptive frequency sampling (AFS) algorithm for transient BEM analysis, in which the sampling frequency is not equally-spaced but adaptively determined according to the characteristics of the computed frequency-domain responses. This work is inspired by the fact that the frequency spectrums of most real structures are often smooth apart from certain peaks; and consequently, the equally-spaced sampling frequencies used in the FSM is not optimal in terms of computational efficiency. The adaptive algorithm is built upon the fitting of the frequency response functions using rational functions, which can be efficiently solved by using the vector fitting method proposed in \cite{GS99}; see Sections \ref{S-fdvf} and \ref{S-adaptive} for the descriptions of the frequency-domain rational fitting and the adaptive algorithm, respectively.

To the knowledge of the present authors, the present adaptive algorithm is unique in solving transient elastodynamic problems using the frequency-domain approach. It assembles the time-domain responses by using the linear combinations of the frequency-domain responses at some ``typical'' frequencies; in this sense, it bears a similarity to the method of modal superposition in finite element analysis. The difference is that in the finite element method these ``typical'' frequencies are the natural frequencies obtained by eigenvalue analysis. However, in BEM we can not afford to solve large-scale nonlinear eigenvalue problems for the natural frequencies even to date \cite{Beyn12}. Our adaptive algorithm works by extracting some ``typical'' frequencies based on the computed frequency response functions.

Finally, we mention that the idea of reducing the number of solutions in frequency domain in preparation for the numerical inversion is not new. It was first employed by Roesset and Kausel in conjunction with the Fourier transform in \cite{kausel1975dynamic} and later by Beskos et al in conjunction with the Laplace transform in \cite{beskos1984dynamic}. However, in all those works equally-spaced sampling frequencies were used and the number of frequency-domain solutions was reduced by using suitable value of the damping parameter $\eta$. This is similar to the FSM. From this point of view, our work builds upon the method used in, e.g., \cite{beskos1984dynamic}. We achieve further reduction of the frequency-domain solutions by using more efficient, unequally-spaced sampling frequencies, which are adaptively determined according to the actually responses.

Our numerical experiments indicate that, by using the proposed adaptive algorithm instead of the FSM, the number of sampling frequencies can be reduced by a factor of 2 or more,
while still preserving the accuracy of the time-domain solutions; see Section \ref{S-ne}. This implies considerable saving of the computational time in large-scale transient elastodynamic BEM analysis.

\section{Basic theories and methods}\label{S-fbem-fda}

\subsection{Fast BEM for frequency-domain elastodynamics}\label{S-fbem}

Let $\Omega \in \R^3$ denote the region of space occupied by a three-dimensional
elastic solid with isotropic constitutive properties
defined by Lam\'e constants $\lambda$ and $\mu$, Poisson's ratio $\nu$ and mass density $\rho$.  The speeds of S and P elastic
waves are denoted by
$c_s
= \sqrt{\mu / \rho}$, $c_p = \sqrt{\lambda + 2\mu /
\rho}$.
Assume that the body force vanishes and the initial displacements and velocities are both zero. Then the frequency domain boundary integral equation (BIE) for elastodynamic problem reads
\begin{equation}\label{eq:bie}
  \begin{aligned}
	c_{ij}(\bbx) u_j(\bbx, s) + (\mathrm{P.V.}) \int_{\Gamma} T_{ij} (\bbx,\bby, s) u_j(\bby, s) \diff \Gamma_{\bby}
	= \int_{\Gamma} U_{ij} (\bbx,\bby, s) \sigma_j(\bby, s) \diff \Gamma_{\bby}, \quad \bbx \in \Gamma
  \end{aligned}
\end{equation}
where, $s$ is the complex frequency; $u_j$ and $\sigma_j$ are components of displacements and tractions in the frequency domain, respectively; (P.V.) indicates a Cauchy principal value (CPV) of the singular integral; the free-term $c_{ij}(\bbx)$ is equal to $0.5\delta_{ij}$ for a smooth boundary at $\bbx$; $U_{ij} (\bbx,\bby; s)$ and $T_{ij} (\bbx,\bby; s)$ denote the displacement and
traction fundamental solutions which can be found in many text books and thus are omitted here.

In this paper, the frequency-domain BIE \eqref{eq:bie} is solved by using the locally-corrected Nystr\"om BEM \cite{COS98} based on curved quadratic elements.
In the numerical implementation, the boundary $\Gamma$ is partitioned into $n_e$ curved triangular quadratic elements.
The 6-point Gauss quadrature rule on triangle is used in evaluating regular element integrals. Thus, the nodes of the Nystr\"om method on each element are the points of the 6-point Gauss rule, and there are totally $n_k=6\cdot n_e$ in the boundary element mesh.
By adopting the Nystr\"om discretization to all the boundary integrals associated with the components of the kernels $U_{ij}$ and $T_{ij}$ and enforcing the boundary conditions, one can finally obtain a linear system of equations
\begin{equation}\label{eq:Aa_b}
  \bbA(s) \bba(s) = \bbb(s),
\end{equation}
where, the $N$ by $N$ ($N =3n_k$) system matrix $\bbA$ and the $N$-vector $\bbb$ are known, and the $N$-vector $\bba$ collects the unknown nodal displacement and traction components which can be obtained by solving the system. Note both the matrix and vectors in \eqref{eq:Aa_b} are functions of the frequency $s$.
All the nearly singular integrals in the Nystr\"om discretization are evaluated by using the usual recursive subdivision quadrature procedure, and all the weakly and strongly singular integrals are computed by using the method recently proposed in \cite{RWX14}.
The linear system \eqref{eq:Aa_b} is solved iteratively by using the generalized minimal residual method (GMRES).
The evaluation of the matrix-vector product is accelerated by using the kernel-independent fast multipole method (KIFMM) \cite{cao2015kernel}. 

\subsection{Frequency-domain approach for transient analysis}\label{S-prob-fda}

In general, by frequency-domain approach (FDA) we mean a method for computing the time-domain responses via the inverse Laplace transform of their frequency-domain counterparts. To be more specific, let $h(s)$ denote a frequency-domain function, which can be the displacement or traction component in BIE \eqref{eq:bie}. Its time-domain counterpart, denoted by $\hat h(t)$, can be expressed by using the Bromwich contour integral for inverse Laplace transform
\begin{equation}\label{eq:invLap}
\hat h(t) = {1 \over 2 \pi \im} \int^{\eta + \im\infty }_{\eta - \im\infty } h(s) e^{st} \diff s,
\end{equation}
where, $s=\eta + \im \omega$, with $\omega$ being the circular frequency; the abscissa of convergence $\eta$ is a real constant chosen to put the
contour to the right of all singularities in $h(s)$, which will be discussed later.

The numerical computation of the Bromwich contour integral \eqref{eq:invLap} is in general an ill-posed problem. This difficulty has led to the diversity of viable numerical approaches
in the literature; see, e.g., \cite{Cohen07}. In \cite{K13} the performance of five different approaches, including the Gaver-Stehfest method, Schapery method, Weeks method, Talbot method and the Fourier series method (FSM), is compared for BEM applications. In large-scale BEM simulations since the solution of linear system \eqref{eq:Aa_b} is often computationally very expensive, the FSM should be the most economical and robust one in all these five methods.
The primary motivation of this paper stems from the promotion of the computational efficiency of the FSM for solving large-scale problems. However, the proposed method in Section \ref{S-adaptive} is more like a new frequency-domain approach to transient analysis. To inspire the new method, the basic idea of the FSM is briefly summarized.

Let $\Delta \omega$ and $\Delta t$ be the circular frequency and time resolutions,
and $T$ be the time period of the transient response.
Given the number of sampling points $N_{\rm{s}}$, one has the basic relations
$\Delta \omega = {2 \pi \over T}$ and $\Delta t = {T \over N_{\rm{s}}}$, and the maximum frequency is $\Omega_{\rm{s}} = N_{\rm{s}} \Delta \omega /2 = \pi / \Delta t$.
Frequency-domain BEM analyses are performed at the first $(N_{\rm{s}}/2 + 1)$ equally-spaced sampling frequencies $s_k= \eta + \im k \Delta \omega, \,(k=0,\cdots,N_{\rm{s}}/2$) to obtain the frequency-domain responses, represented by $h(s_k)$ in \eqref{eq:invLap}. The responses corresponding to the last $(N_{\rm{s}}/2 - 1)$ frequencies are determined by using the conjugate
symmetric property of the DFT $h(s_k) = h^*(s_{N_{\rm{s}}-k}), \, k={N_{\rm{s}}/ 2}+1, \cdots,  N_{\rm{s}}-1,$
with $h^*$ denoting the complex conjugate of $h$. In the remainder of this paper, we use $N_{\rm{c}}$ to represent the total number of frequencies at which the frequency-domain BEM simulations have to be conducted. Obviously, in the FSM, $N_{\rm{c}} = N_{\rm{s}}/2 + 1$.

The time-domain responses can be computed by truncating the infinite integral in \eqref{eq:invLap} into a finite one and using the trapezoidal rule,
\begin{equation}\label{eq:IDFT-MFT}
	  \hat h(n \Delta t) = e^{\eta n \Delta t} \sum^{N_{\rm{s}} - 1}_{k=0}  W(s_k) h(s_k) e^{ 2\pi \im {n k \over N_{\rm{s}}} }, \quad (n=0,\cdots,N_{\rm{s}} - 1).
\end{equation}
The parameter $\eta$ can be determined by $\eta = \kappa { \ln 10 / T }$,
with $\kappa$ being a user defined constant that controls the trade-off between the effective elimination the aliasing error and the magnification of the truncation error due to the truncation of high-frequency components. We suggest to use $\kappa \in [2.0, 3.0]$. $W(s_k)$ is a frequency-domain window function which is employed to reduce the truncation error and thus the Gibbs oscillation in the time-domain results. In elastodynamic BEM applications, the Hanning windows given by $W(s_k) = 0.5 \left[ 1+ \cos \left( {2 \pi k \over N_{\rm{s}}} \right) \right]$ is a good choice \cite{XYW13}.

We notice that in most situations the equally-spaced sampling frequencies used in the FSM is not optimal in terms of computational accuracy and efficiency, because the frequency responses of the real structures are often not uniformly smooth but contains nonuniformly distributed peaks. To well represent the frequency spectrum, it is more natural to allocate more sampling points around the peaks and less in the smooth area. As such, we expect to keep
the number of frequencies $N_{\rm{c}}$ reasonably small, while still to preserve the accuracy of the time-domain solutions.
This is especially important in large-scale simulations where conducting one BEM analysis by forming and solving linear system \eqref{eq:Aa_b} is computationally very expensive.
In section \ref{S-adaptive} we will propose an adaptive frequency-domain sampling algorithm based on the rational fitting of the frequency responses. By using this algorithm one can sequentially select optimal sampling frequencies according to the computed frequency responses at hand, leading to a reasonably minimized number $N_{\rm{c}}$ without sacrificing the accuracy of the time-domain results when compared with the FSM.

\subsection{Frequency domain rational fitting}\label{S-fdvf}

It is a common knowledge that the frequency response function of a dynamic system can be well approximated by rational functions. A number of algorithms have been proposed for this purpose, see e.g., the Stoer-Bulirsch algorithm \cite{SB80}, asymptotic wave form evaluation \cite{PR90}, Pad\'{e} approximations \cite{HTF12} and the vector fitting (VF) algorithm \cite{GS99}. However, many of those algorithms are not suitable for large-scale BEM analysis. For example, the Pad\'{e} interpolating algorithm \cite{HTF12} involves high order derivatives of the system matrix, of which both the computational cost and memory requirement are huge in large-scale BEM. The VF is a robust algorithm for solving rational fitting problems. It permits to identify frequency spectrums directly from computed frequency responses. The adaptive frequency sampling algorithm in this paper will be based on the VF algorithm.

Below, vector function $\bbf (s) = \{f_1(s), f_2(s), \cdots, f_K(s) \}^T$ will be used to denote a subset of nodal responses in the solution vector $\bba (s)$ of the BEM system \eqref{eq:Aa_b}, which will be approximated by using the VF. The total number of functions in $\bbf (s)$, $K$, can be a small number that will be specified below. Function $f_k(s)$ can be the displacement or traction response at certain point in certain direction. For real structures, it is reasonable to approximate $\bbf (s)$ by using partial fraction expansion
\begin{equation} \label{vf-f}
\bbf (s) = \sum^M_{m=1} { \bbr_m \over s-a_m }
\end{equation}
based on simulated frequency responses $\{s_j, \bbf(s_j)\}_{j=1}^J$.
In \eqref{vf-f} $\{a_m\}$ is the set of common poles for all the response functions in $\bbf$, $\bbr_m$ are vectors of the residues, $M$ is the order of the approximation.

The VF proposed
by Gustavsen and Semlyen \cite{GS99} is an accurate and stable algorithm for solving the rational fitting problem \eqref{vf-f}. It consists of two basic steps: \emph{pole identification} and \emph{residue identification}. In the pole identification step, the poles $\{a_m\}$ are extracted from the frequency response data and a given set of initial poles $\{\bar a_m \}$ iteratively.
This procedure is known to have good convergence properties (usually
converges in a few iterations) \cite{G06}.
After the poles $\{a_m\}$ have been identified, the residues $r_{k,m}$ can be easily calculated by solving the corresponding
least-squares problem with known poles. To make sure that the
transfer function has real-valued coefficients, the linear problem
is constructed in such a way that a complex-conjugate pole
pair results in a complex-conjugate residue pair.

There are two issues that have substantial influence on the accuracy of the algorithm: the initial poles $\{\bar a_m\}$ and the treatment of unstable poles occurred during iterations. In the present paper the initial poles are always chosen as equally spaced in the interested frequency band. There are two ways to handle the unstable poles. One is to force the unstable poles stable by flipping them into the left half-plane (inverting the sign of their real parts). Another way is simply kept the unstable poles unchanged. Our frequency-domain fitting for elastodynamics shows that the first approach could often result in large errors while the second always performs well. We therefore retain the unstable poles in this paper.


Here we notice that in this paper we only select a small number $K$ of response functions in the pole-identification process. However, the poles obtained by these $K$ functions are used to fitting the responses of all the DOFs in the BEM modal. The rationality of this will be discussed in Section \ref{SS-ORF} and numerically verified in Section \ref{S-ne}. Given the poles $\{a_m\}$, let $\{r_m\}$ be the residues of frequency response function $\{h(s_j)\}_{j=1}^J$, then the time-domain response of $h(s)$ can be obtained by the inverse Laplace transform of the rational expansion \eqref{vf-f} which has close form
\begin{equation}\label{eq:vf-td}
\hat h (t) = \sum^M_{m=1} r_m e^{-a_m t }.
\end{equation}
This is a continuous time-domain response function, while by using the FSM one can only obtain $\hat h (t_n)$ at the equally-spaced time steps \eqref{eq:IDFT-MFT}.

\section{Adaptive frequency sampling}\label{S-adaptive}

The adaptive frequency sampling (AFS) algorithm proposed in this paper is an iterative process. It consists of three stages. The first stage is initialization in which one should first select a set of observation response functions (ORFs) $\bbf(s)$ from the solution vector $\bba(s)$ of BEM system \eqref{eq:Aa_b} for vector fitting. A series of initial frequencies $\bbs^{(0)}$ and the corresponding values of $\bbf(s)$ at $\bbs^{(0)}$, denoted by $\bbf(\bbs^{(0)})$, are also needed to start the AFS. Notice that $\bbf(\bbs^{(0)})$ is obtained by performing BEM analysis at the frequencies in $\bbs_0$. In the second stage new frequencies are iteratively selected. Given the current frequencies $\bbs^{(i)}$ and responses $\bbf(\bbs^{(i)})$ at step $i$, two vector fittings with different orders are conducted, and the new sampling frequency $s^{(i)}_{\rm{new}}$ is determined according to the discrepancy of the two fittings. Once $s^{(i)}_{\rm{new}}$ is obtained, a new BEM analysis is performed and the response data $\bbf(s^{(i)}_{\rm{new}})$ are collected. The third stage is convergence examination. Some types of errors should be introduced and checked in each AFS step to see if the AFS process can be terminated. In the following subsections, technique details of each stage will be described.

\subsection{Selection of observation response functions} \label{SS-ORF}

First and foremost, one needs to fix on how many functions should be included in $\bbf(s)$ and what are they. A basic principle is that $\bbf(s)$ and its discretization $\bbf(\bbs)$ at sampling frequencies $\bbs$ should reflect the main characteristics of the structure in the given frequency band; or more specifically, the poles $\{a_m\}$ and the rational basis function $ 1 / (s-a_m) $, which are extracted from $\bbf(\bbs)$ by the VF algorithm, should well fit the frequency spectrum of the structure in entire frequency interval $[0,\, \Omega_{\rm{s}}]$. By assuming that the poles obtained by the VF algorithm are sufficiently accurate, it is natural to expect that more dynamic information of the structure will be reflected in $\{a_m\}$ if more response functions are included in $\bbf(s)$. However, one should note that the computational cost of the pole-identification process in VF increases quadratically with the number of functions in $\bbf(s)$. Is it possible to obtain good poles $\{a_m\}$ using only a small number of ORFs? The answer seems to be positive. In our numerical experiments, we randomly selected $K=5\sim 10$ ORFs from the BEM model.
We found that the poles $\{a_m\}$ obtained by these ORFs can always fit all the nodal frequency responses very well, provided that the number of sampling frequencies are large enough. Besides the random selection, one can also choose the most interested response functions as ORFs.

A rigorous justification of the above observation is difficult so far; however, we can draw some hints from the theory of the nonlinear eigenvalue problems. Our explanation is based on the Keldysh's theorem of the nonlinear eigenvalues; see \cite{Beyn12} and the references therein. Indeed, the poles of the frequency-domain responses are mainly determined by the eigenvalues, denoted by $\{\lambda_m\}$, of the system matrix $\bbA(s)$ in \eqref{eq:Aa_b}. Suppose that $\bbA(s)$ has only simple eigenvalues, then it follows from the Keldysh's theorem that $\bbA(s)^{-1}$ can be expressed as a linear combination of the rational functions with poles $\{\lambda_m\}$,
$[\bbA(s)]^{-1} = \sum_{m} {S_{m} \over (s-\lambda_m)}$,
with $S_{m}$ being the coefficient matrices of the same size as $\bbA(s)$.
By further invoking the solution of linear system \eqref{eq:Aa_b} $ \bba(s) = \bbA(s)^{-1}\bbb(s)$, and assuming that $\bbb(s)$ is an analytic vector function, then it is reasonable to conclude that there is a high probability for each element $f_k(s)$ of $\bba(s)$ to take almost all the information of the eigenvalues of $\bbA(s)$.
As such, one can expect to identify these common poles of $\bba(s)$ by randomly selecting several functions from $\bba(s)$. Our numerical experiments indicate that, by using the VF algorithm, we can indeed identify some eigenvalues of matrix $\bbA(s)$. However, there are also many poles of VF that are not the eigenvalues of $\bbA(s)$ but necessary to keep the accuracy of the rational fitting of the frequency spectrum in $[0,\, \Omega_{\rm{s}}]$. The explanation of this phenomenon is an interesting topic for further investigations.

\subsection{Determination of new sampling frequencies}\label{S-S-snew}

Notice that most of the notations in this subsection, except the two constants $\alpha_{\rm{H}}$ and $\alpha_{\rm{L}}$ in \eqref{eq:fit-ords}, are associated with certain step of the AFS process; however, since the operations are common in each step, we neglect the step index superscript ``$(i)$'' in the notations.

\begin{figure}
\centering
\includegraphics[width = 0.45\textwidth]{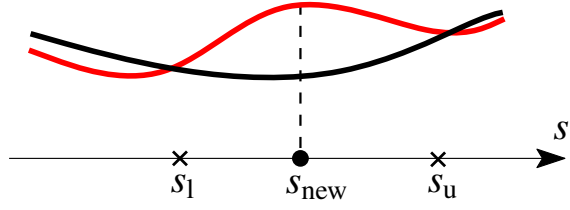}\\
\caption{Strategy for determining the new sampling frequencies. The black and red curves denote the two fittings of low and high orders, respectively. $s_{\rm{l}}$ and $s_{\rm{u}}$ represent the fixed frequencies used for fitting. }
  \label{fig:fitting}
\end{figure}

The central to the AFS is, for a given set of sampling frequencies $\bbs$ and the corresponding function values $\bbf(\bbs)$, how to determine a new frequency $s_{\rm{new}}$ at which a new BEM system $\bbA(s_{\rm{new}}) \bbf(s_{\rm{new}}) = \bbb(s_{\rm{new}})$ will be formulated and solved. Our idea is to approximate data $\left[\bbs,\bbf(\bbs)\right]$ with two rational fittings of different orders $M_{\rm{H}}$ and $M_{\rm{L}}(<M_{\rm{H}})$. Since there is no prior knowledge of the elastodynamic system, we assume that the fitting with higher order $M_{\rm{H}}$ to be a more accurate model of the system. This becomes reasonable as the number of samplings grows. The second fitting with lower order $M_{\rm{L}}$ is used to determine the location of the new frequency $s_{\rm{new}}$. Figure \ref{fig:fitting}(a) is a sketch of this idea in the scalar case ($K=1$). A new frequency $s_{\rm{new}}$ is inserted where the largest discrepancy of the two fittings occurs, i.e.,
\begin{equation}\label{eq:snew}
s_{\rm{new}} = \underset{s \in [0, \Omega_{s}] }{ \mathrm{argmax} } \, |f_{\rm{H}}(s) - f_{\rm{L}}(s) |,
\end{equation}
where, $f_{\rm{H}}(s)$ and $f_{\rm{L}}(s)$ denote the two fittings of $f(s)$.

The determination of the new sampling frequency $s_{\rm{new}}$ in the vector case $K>1$ is more complicated, because the frequencies corresponding to the largest discrepancies of the $K$ functions in $\bbf(s)$ are often not coincident with each other. Thus, one has to decide which frequency should be chosen as $s_{\rm{new}}$. The solution is not unique. Our idea is to determine $s_{\rm{new}}$ by the response function in $\bbf(s)$ which has the largest relative error. To be concrete, define the relative error of $f_k (s)$ in the current fitting as
\begin{equation}\label{eq:ferr}
 e_k =  { \max\limits_{s \in [0, \Omega_{s}]}|f_{k, \rm{H}}(s) - f_{k, \rm{L}}(s) |  \over \max\limits_{s \in [0, \Omega_{s}] } | f_{k, \rm{H}}(s) | }.
\end{equation}
The index of the function that has the largest error is given by
$$
k_{\max} = \underset{k \in \mathcal {I}_{\rm{o}}}{\mathrm{argmax}} \, e_k,
$$
with $\mathcal {I}_{\rm{o}}$ being the index set of the functions in $\bbf$, which is a subset of DOFs in the BEM analysis. Then, we determine $s_{\rm{new}}$ by using the two fittings of function $f_{k_{\max}}(s)$ in the same way as the scalar case; that is,
\begin{equation}\label{eq:snew1}
s_{\rm{new}} = \underset{s \in [0, \Omega_{s}] }{ \mathrm{argmax} } \, |f_{k_{\max}, \rm{H}}(s) - f_{k_{\max}, \rm{L}}(s) |.
\end{equation}

The fitting orders $M_{\rm{H}}$ and $M_{\rm{L}}$ are determined according to the number of sampling frequencies $J$ in the current AFS step as
\begin{equation}\label{eq:fit-ords}
M_{\rm{H}} = \left\lfloor {J \over \alpha_{\rm{H}} } \right\rfloor, \quad M_{\rm{L}} = \left\lfloor {J \over \alpha_{\rm{L}} } \right\rfloor,
\end{equation}
where, $\alpha_{\rm{H}}$ and $\alpha_{\rm{L}}$ are two constants. In this paper we always use around 16 initial frequencies in the AFS, and we find that exellent performance can be attained by setting $\alpha_{\rm{H}}=2.0$ and $\alpha_{\rm{L}}=2.3$.

Once the new frequency $s_{\rm{new}}$ is obtained, a BEM analysis at $s_{\rm{new}}$ is performed and the values of the ORFs, denoted by $\bbf(s_{\rm{new}})$, can be obtained. Then, the set of sampling frequencies is updated as
$\bbs \leftarrow \bbs \cup {s_{\rm{new}}}$,
which will be used in the next step of the AFS.

\subsection{Time-domain response and stopping criterion of AFS}

Given the data $\left[\bbs^{(i)},\bbf(\bbs^{(i)})\right]$ at the $i$-th step, the time-domain responses of the ORFs can be computed as
\begin{equation}\label{eq:vf-td-afs}
\hat \bbf^{(i)} (t) = \sum^{M_{\rm{H}}}_{m=1} \bbr_m e^{-a_m t },
\end{equation}
with $\{\bbr_m\}$ and $\{a_m\}$ being obtained the order-$M_{\rm{H}}$ fitting of the data. To compute the time-domain response of any other component of the solution vector $\bba$ in \eqref{eq:Aa_b}, one should first perform a residue identification operation to obtain $\{r_m\}$ using the common poles $\{a_m\}$; see Section \ref{S-fdvf}.

The convergence of the AFS process is controlled by certain error tolerances. A pertinent definition of the convergence error is crucial since it can avoid the possibilities
of under- or over-sampling. Two errors are defined and used in this paper. The first one, denoted by $E_1$, is the error of the time-domain solutions in two successive AFS iterations. $E_1$ can be estimated by the time-domain responses of a subset of functions in the solution vector $\bba$. Let $\mathcal {I}_{\rm{t}}$ denote the index set of these functions and let $\bar f^{(i)}_k (t), (k \in \mathcal {I}_{\rm{t}})$ represent their time-domain responses computed by \eqref{eq:vf-td-afs} in the $i$-th step. Then, $E_1$ is defined as
\begin{equation}\label{eq:stop-E1}
E^{(i)}_1 = \max\limits_{k \in \mathcal {I}_{\rm{t}} } { \int_0^T [ \bar f^{(i)}_k (t) - \bar f^{(i-1)}_k (t) ]^2 \diff t  \over  \int_0^T [ \bar f^{(i)}_k (t) ]^2 \diff t  }.
\end{equation}
The second error $E_2$ is the fitting error in the frequency domain. It can be estimated by using the relative error $e_k$ defined in \eqref{eq:ferr}; that is,
\begin{equation}\label{eq:stop-E2}
E^{(i)}_2 = \max \{e^{(i)}_k, \, k \in \mathcal {I}_{\rm{t}} \}.
\end{equation}
Notice that in \eqref{eq:ferr} the superscript of $e^{(i)}_k$ is neglected. Also notice that in \eqref{eq:stop-E2} $E_2$ is estimated using the same set of DOFs $\mathcal {I}_{\rm{t}}$ as $E_1$; this is, in general, not a necessity. For DOFs which are not in the set $\mathcal {I}_{\rm{o}}$ of the ORFs, the corresponding frequency-domain fittings can be obtained by using the common poles extracted by the ORFs and performing a residual identification step.

By using the above errors, the convergence of the AFS process can be defined and checked as follows. If the AFS process is converged at step $i$, then
the two errors $E^{(i')}_1$ and $E^{(i')}_2$ in all the following steps $(i'=i,i+1,i+2,\cdots)$ must be less than some user-defined thresholds $\bar E_1$ and $\bar E_2$, respectively; that is,
\begin{equation}\label{eq:converge}
E^{(i')}_1 < \bar E_1 \quad \mathrm{and} \quad E^{(i')}_2 < \bar E_2, \quad i'=i,i+1,i+2,\cdots.
\end{equation}
In implementation, the two errors are monitored in each step. The AFS process can be terminated once condition $E^{(i)}_1 < \bar E_1$ and $E^{(i)}_2 < \bar E_2$ is satisfied at certain step $i$.
However, this can not guarantee the satisfaction of the condition for all the following steps $(i+1, i+2, \cdots)$, and thus may lead to possible under-sampling.
To reduce the possibility of under-sampling, one can keep the AFS process going for a small number, denoted by $J_{\rm{c}}$, of steps after step $i$ and terminate the AFS process only if the condition can be satisfied in all the following $J_{\rm{c}}$ steps. In our implementation, $J_{\rm{c}}=3$ is used.

\subsection{Implementation details of the AFS algorithm} \label{S-S-imp}

The AFS algorithm is summarized as follows.

\textit{
Algorithm (AFS)
\begin{enumerate}[(1)]
  \item Initialization
  \begin{itemize}
    \item Determine the set $\mathcal {I}_{\rm{o}}$ of ORFs $\bbf(s)$ for VF and the set $\mathcal {I}_{\rm{t}}$ of functions for testing the convergence.
    \item Set values for $\bbs^{(0)}$, $\alpha_{\rm{H}}$ and $\alpha_{\rm{L}}$ in \eqref{eq:fit-ords}, $\bar E_1$ and $\bar E_2$ in \eqref{eq:converge}, $J_{\rm{c}}$, etc.
    \item For each of the $J_0$ frequencies in $\bbs^{(0)}$, perform BEM analysis and obtain $\bbf(\bbs^{(0)})$.
    \item Compute $E^{(0)}_1$ and $E^{(0)}_2$, and set step index $i=0$.
  \end{itemize}
  \item Adaptive frequency sampling process\\
  \indent While $E^{(i)}_1 > \bar E_1$ or $E^{(i)}_2 > \bar E_2$, do
  \begin{itemize}[---]
    \item Find $s^{(i)}_{\rm{new}}$ by using data $\left[\bbs^{(i)},\bbf(\bbs^{(i)})\right]$; see Section \ref{S-S-snew}.
    \item Perform BEM analysis at $s^{(i)}_{\rm{new}}$ and obtain the values of ORFs $\bbf(s^{(i)}_{\rm{new}})$.
    \item Set $\bbs^{(i+1)} \leftarrow \bbs^{(i)} \cup {s^{(i)}_{\rm{new}}}$.
    \item Compute $E^{(i+1)}_1$ and $E^{(i+1)}_2$, and set $i \leftarrow i+1$.
  \end{itemize}
    \indent End while
  \item Convergence validation\\
  \indent For $i'=i+1, \cdots, i+J_{\rm{c}}$, do
    \begin{itemize}[---]
    \item Perform adaptive sampling as in Stage (2), and obtain $\bbs^{(i')}$ and $\bbf(\bbs^{(i')})$.
    \item Compute $E^{(i')}_1$ and $E^{(i')}_2$.
    \item if $E^{(i')}_1 > \bar E_1$ or $E^{(i')}_2 > \bar E_2$, then set $i=i'$, go to stage (2).
  \end{itemize}
  \indent End for
  \item Post-processing\\
  \indent For example, computing the time-domain responses via \eqref{eq:vf-td}.
\end{enumerate}
}

In this paper, the initial sampling frequencies $\bbs^{(0)}$ are set to be the Chebyshev points in $[0\; \Omega_{\rm{s}}]$, i.e.,
$$
s_j = \eta + \im \Omega_{\rm{s}} \left[ 1-\cos {j \over 2(J_0-1)} \right], \quad j = 0,1, \cdots, J_0-1,
$$
with $J_0=16$. The other parameters are set as: $\alpha_{\rm{H}}=2.0$, $\alpha_{\rm{L}}=2.3$, $\bar E_1=10^{-4}$, $\bar E_2=2\cdot 10^{-3}$, $J_{\rm{c}}=3$. These values were found to be good in achieving the smallest number of frequencies $N_{\rm{c}}$.

\section{Numerical examples}\label{S-ne}

The performance of the AFS algorithm proposed in this paper is assessed here for three
representative examples. The first one is a prismatic rod subject to a step traction pulse, which is a benchmark in testing the capability of time-domain solvers. The second one is a thick plate with a hole in middle subject to an impact loading, for which experimental results are available for reference. In the third example elastic wave propagation in a porous solid model is simulated, in order to manifest the peculiarity of the algorithm for models involving complicated geometries and large number of DOFs.

The main advantage of the proposed AFS algorithm, in contrary to the FSM with equally-spaced frequencies, is the effective reduction of $N_{\rm{c}}$, while still keeping the same precision in the time-domain results. To demonstrate this feature, the problems in the examples are solved
by using both AFS and FSM, and the respective $N_{\rm{c}}$ of these two methods are compared under the condition that their time-domain results are almost the same.

All the simulations are conducted by using the fast Nystr\"om BEM in Section \ref{S-fbem} on a computer with a Xeon 3.0GHz CPU and 30GB RAM.
The linear systems are solved by using GMRES in conjunction with the block diagonal preconditioner. The damping coefficient is set as $\kappa = 3.0$, and the Hanning window is used in computing the time-domain solution.

\subsection{Example 1: Prismatic rod subject to a step loading}\label{S-S-NE1}

\begin{figure}[hbt]
\centering
\begin{minipage}[b]{.5\textwidth}
\centering
\epsfig{figure=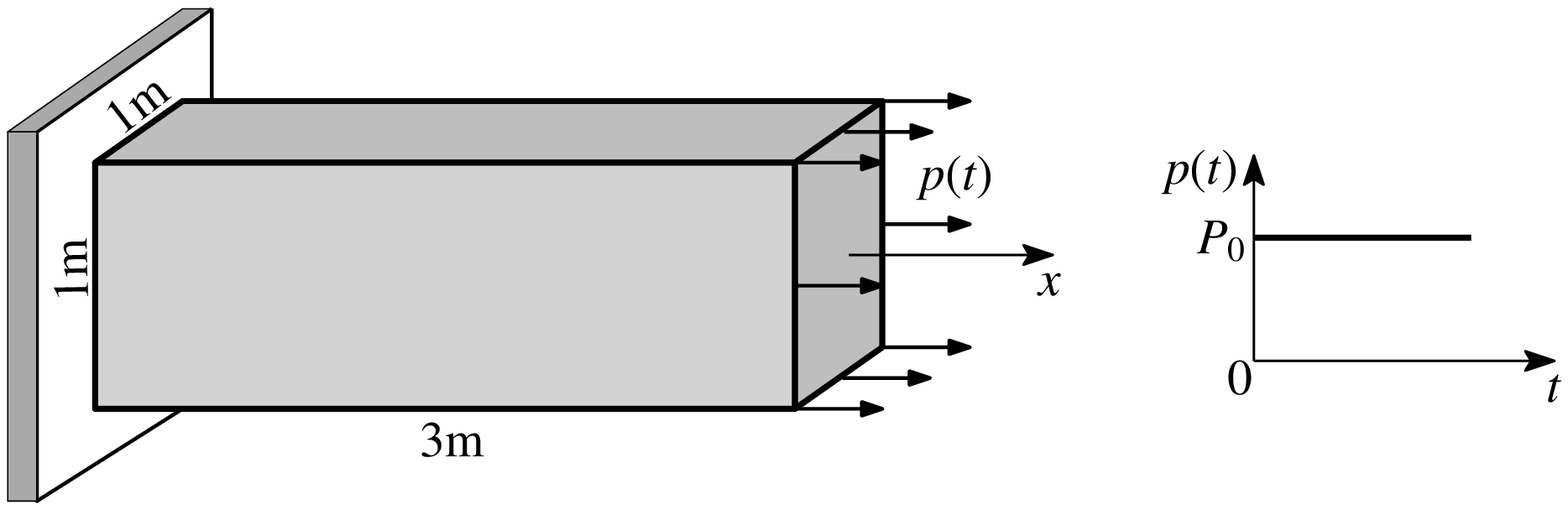,width=0.9\textwidth}\\
\end{minipage}%
\begin{minipage}[b]{.5\textwidth}
\centering
\epsfig{figure=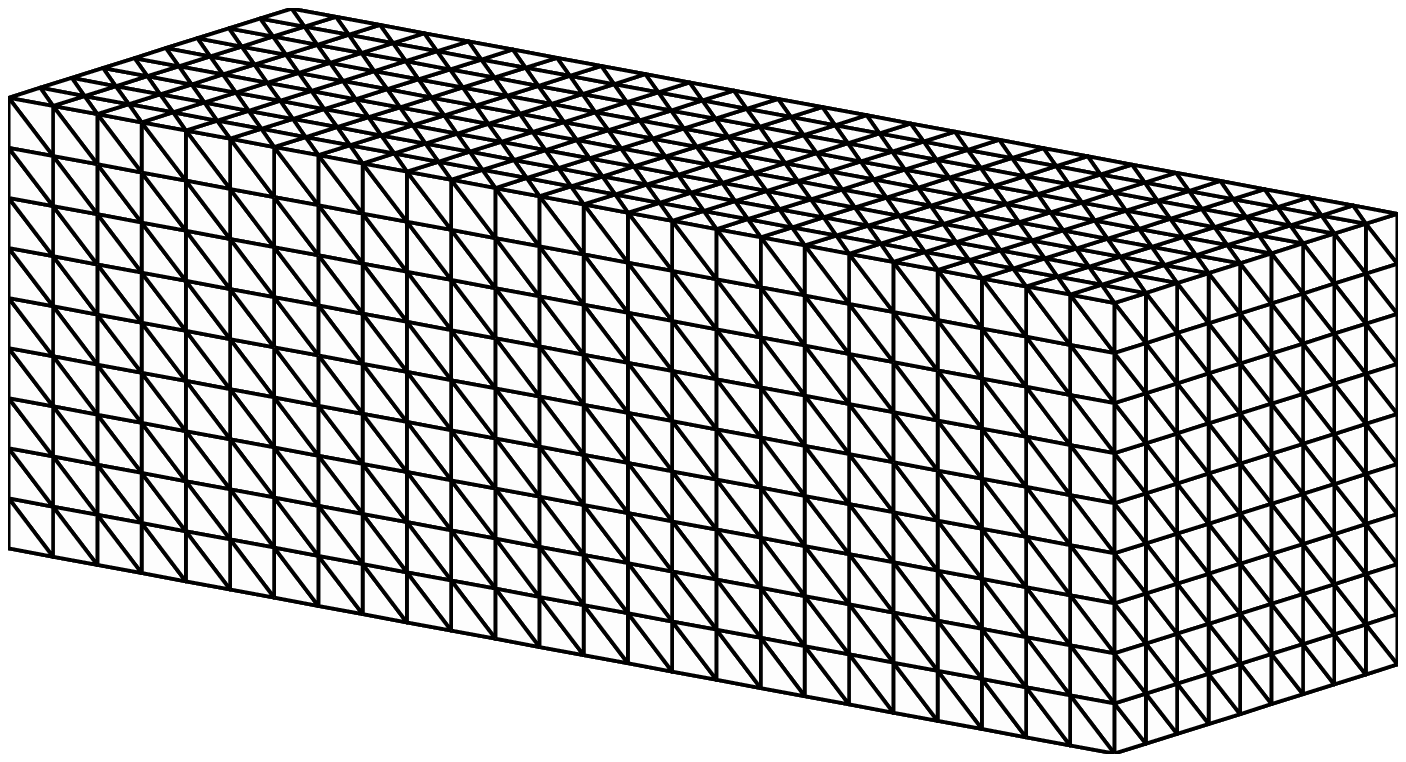,width=0.9\textwidth}\\
\end{minipage}
\caption{Prismatic rod subject to a step traction pulse. Right: model and size; left: mesh.}
\label{fig-rod}
\end{figure}

As the first test case, a prismatic rod with dimensions of $3\rm{m}\times 1\rm{m} \times 1\rm{m}$ as sketched in Figure \ref{fig-rod} is simulated. The rod is fixed at the left end and subjected to a Heaviside traction $p(t)$ = $P_0 H(t)$ at the right end, where $P_0 = -10^6$ N/m$^2$ and $H(t)$ is the Heaviside function. By setting the Poisson's ratio of the material to zero, the 3D problem is reduced to a 1D problem, of which the analytical solution is available, see e.g. \cite{XYCZ12}. The Young's modulus $E$ and the density $\rho$ of
the material are set to be $E = 2.11\cdot 10^{11}$Pa and
$\rho = 7.85 \cdot 10^{3}$kg/m$^3$, respectively.

A boundary element mesh with 2000 triangular quadratic elements is used; see Figure \ref{fig-rod}.
The time-domain simulation is performed by using FSM and AFS. The solution period is set as $T=0.0155$s. The controlling parameters of the AFS are set according to the description in Section \ref{S-S-imp}. In the FSM, one has $\Delta t = \pi/\Omega_{\rm{s}} = 7.11\cdot 10^{-5}$s, $N_{\rm{s}} = \lceil T/\Delta t \rceil = 218$ and $N_{\rm{c}} = 110$. The frequency-domain BEM analysis is preformed at frequencies $s_k= \eta + \im k \Delta \omega$, $(k=0,\cdots,N_{\rm{s}}/2$), and then the time-domain response is obtained by using \eqref{eq:IDFT-MFT}.
%

Let us first check the precision of the frequency-domain VF, which has essential influence on the efficiency of the AFS. In the FSM simulation, the frequency-domain solutions at the $110$ frequencies have been solved. Here we select 5 points on the surface of the bar and fit the computed displacement and/or traction responses at these points by using the VF algorithm. Three different numbers of poles $(M=108,\;30,\;16)$ are used in VF. In order to show that the poles extracted from the 5 ORFs can also
be used to well approximate the responses of the remaining DOFs on the model,
we randomly select 2 testing points on the surface.
The quality of the fitting is measured by the following error,
$$
E_{\rm{vf}} = \sqrt{ {  \sum_{j=0}^{N_{\rm{s}}/2} [f(s_j) - f_j]^2 \over  \sum_{j=0}^{N_{\rm{s}}/2} f_j^2 } },
$$
with $f(s)$ denoting the rational approximation obtained by the VF algorithm and $f_j$ the fitted responses at the equally-spaced sampling frequencies.

The coordinates of the 5 observation points and 2 testing points, the types of responses and the fitting error $E_{\rm{vf}}$ are listed in Table \ref{tab-NE1-vferr}. It is seen that an accuracy up to 4 significant digits can be achieved by using only 30 poles. Further increasing the number of poles does not evidently improve the quality of the fitting; but decreasing the number of poles makes the fitting dramatically deteriorated. The former is due to the fact that the responses at the 110 frequencies suffer from numerical errors and thus not necessarily coincide with rational functions. The latter is more natural since the order of the expansion is reduced. One can conclude that highly accurate approximations of the frequency-domain responses can be achieved by using the VF algorithm with sufficient number $M$ of poles (i.e., order of expansion in \eqref{vf-f}). Furthermore, our experience indicates that $M$ can always be of moderate value that is much less than the number needed in the equally-spaced sampling. Indeed, this is the primary motivation for proposing the AFS algorithm in this paper.

More interestingly, the above results indicate that, in fact, only around 30 frequency-domain simulations are needed to obtain good time-domain solutions. However, the problem is that one does not know the placement of these optimal frequencies beforehand. Our resolution is to iteratively and adaptively find ``nearly optimal'' frequencies using the AFS algorithm.

\begin{table}[hbt]
\begin{center}
\caption{Errors of the frequency-domain VF. Test---testing response; disp---displacement; trac---traction.}\label{tab-NE1-vferr}
\begin{tabular*}{0.8\textwidth}{@{\extracolsep{\fill}}cllllll@{}}\toprule
Index   & Coordinates  & Type     & Response   & $E_{\rm{vf}},M108$ & $E_{\rm{vf}},M30$ & $E_{\rm{vf}},M16$\\
\hline
1       & $(0.379, 0.187, 3)$           & ORF    &  Disp  & 6.37e-5  & 1.91e-4 & 3.06e-1\\
2       & $(0.409, 0.0453, 0)$          & ORF    &  Trac  & 6.04e-5  & 2.10e-4 & 4.20e-1 \\
3       & $(-0.177, -0.369, 3)$         & ORF    &  Disp  & 8.00e-5  & 2.07e-4 & 3.06e-1  \\
4       & $(-0.146, -0.0657, 0)$        & ORF    &  Trac  & 5.66e-5  & 2.64e-4 & 4.19e-1\\
5       & $(-0.0657, 0.5, 1.53)$        & ORF    &  Disp  & 4.80e-5  & 1.89e-4 & 3.88e-1 \\
6       & $(0.298, -0.177, 0)$          & Test   &  Trac  & 5.25e-5  & 2.18e-4 & 4.19e-1  \\
7       & $(-0.5, 0.156, 1.898)$        & Test   &  Disp  & 5.85e-5  & 2.15e-4 & 3.74e-1\\
\bottomrule
\end{tabular*}
\end{center}
\end{table}


Now we solve the transient response by using the AFS. The first 5 responses in Table \ref{tab-NE1-vferr} are selected as the ORFs in the VF process; all the 7 responses are employed to evaluate the time- and frequency-domain errors $E_1$ \eqref{eq:stop-E1} and $E_2$ \eqref{eq:stop-E2} and check the convergence. The values of the controlling parameters are set to be the suggested values in Section \ref{S-S-imp}. Once convergence is reached, the time-domain solution is computed by \eqref{eq:vf-td} with $M=M_{\rm{H}}$. We compare the AFS with the FSM in terms of the quality of the frequency-domain and time-domain solutions and the total number $N_{\rm{c}}$ of sample frequencies required.

\begin{figure}[hbt]
\centering
\epsfig{figure=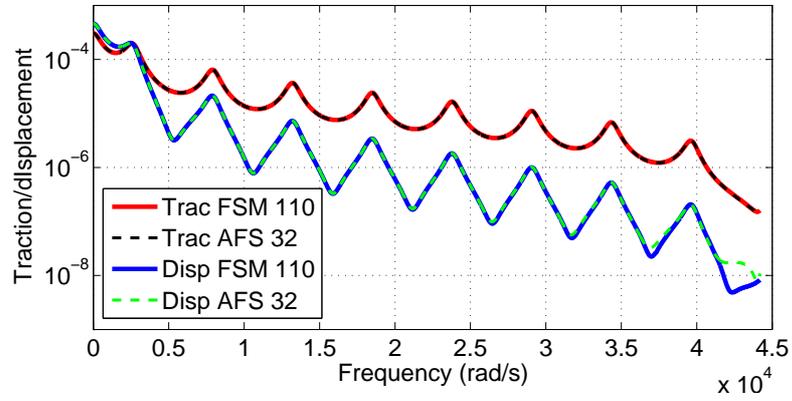,width=0.7\textwidth}
\caption{The amplitudes of the frequency-domain responses of observation points $1$ and $2$.
The traction response is properly amplified to fit the figure. }
\label{fig-n1-fd-obs12}
\end{figure}

\begin{figure}[hbt]
\centering
\epsfig{figure=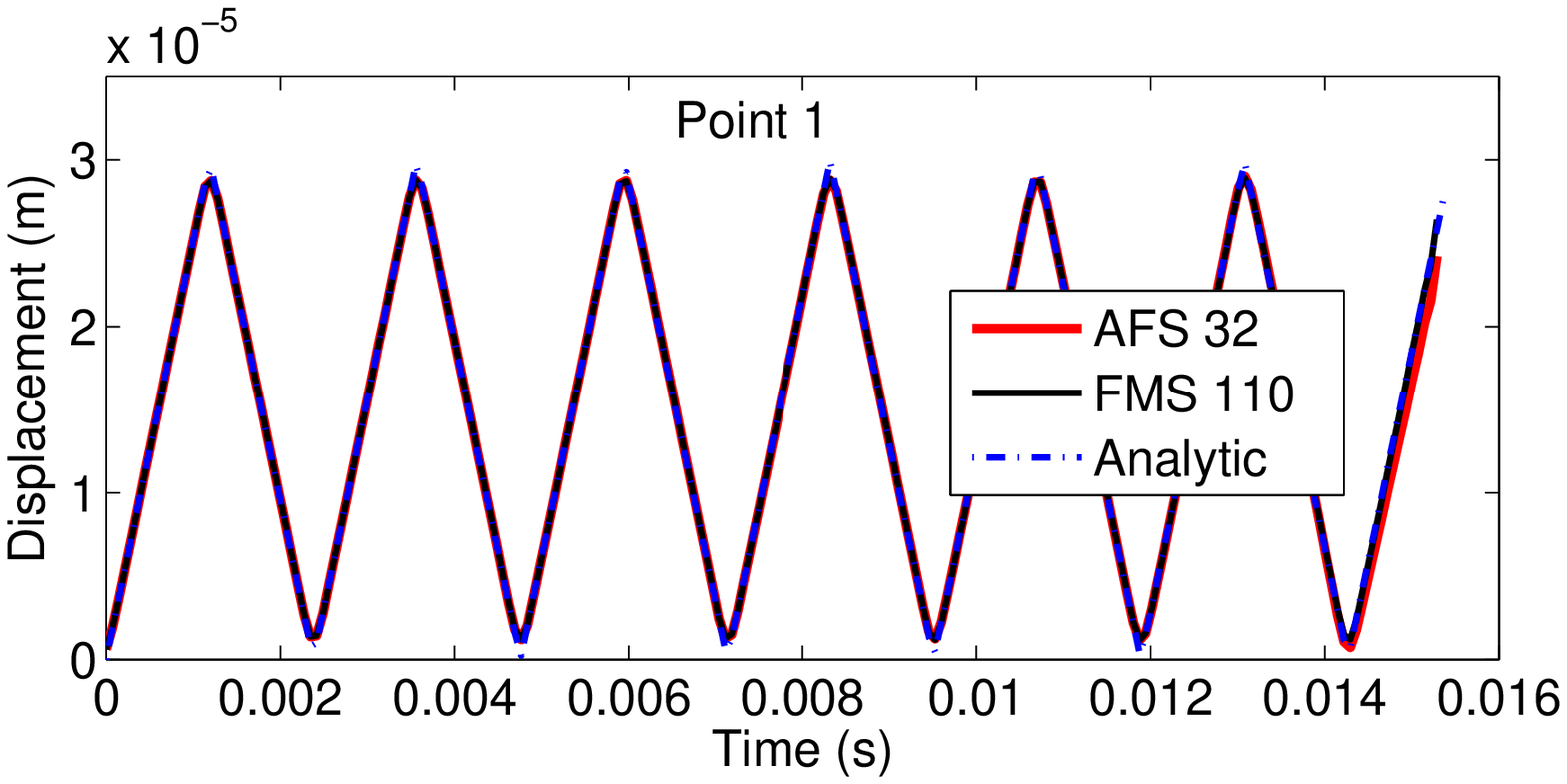,width=0.7\textwidth}
\epsfig{figure=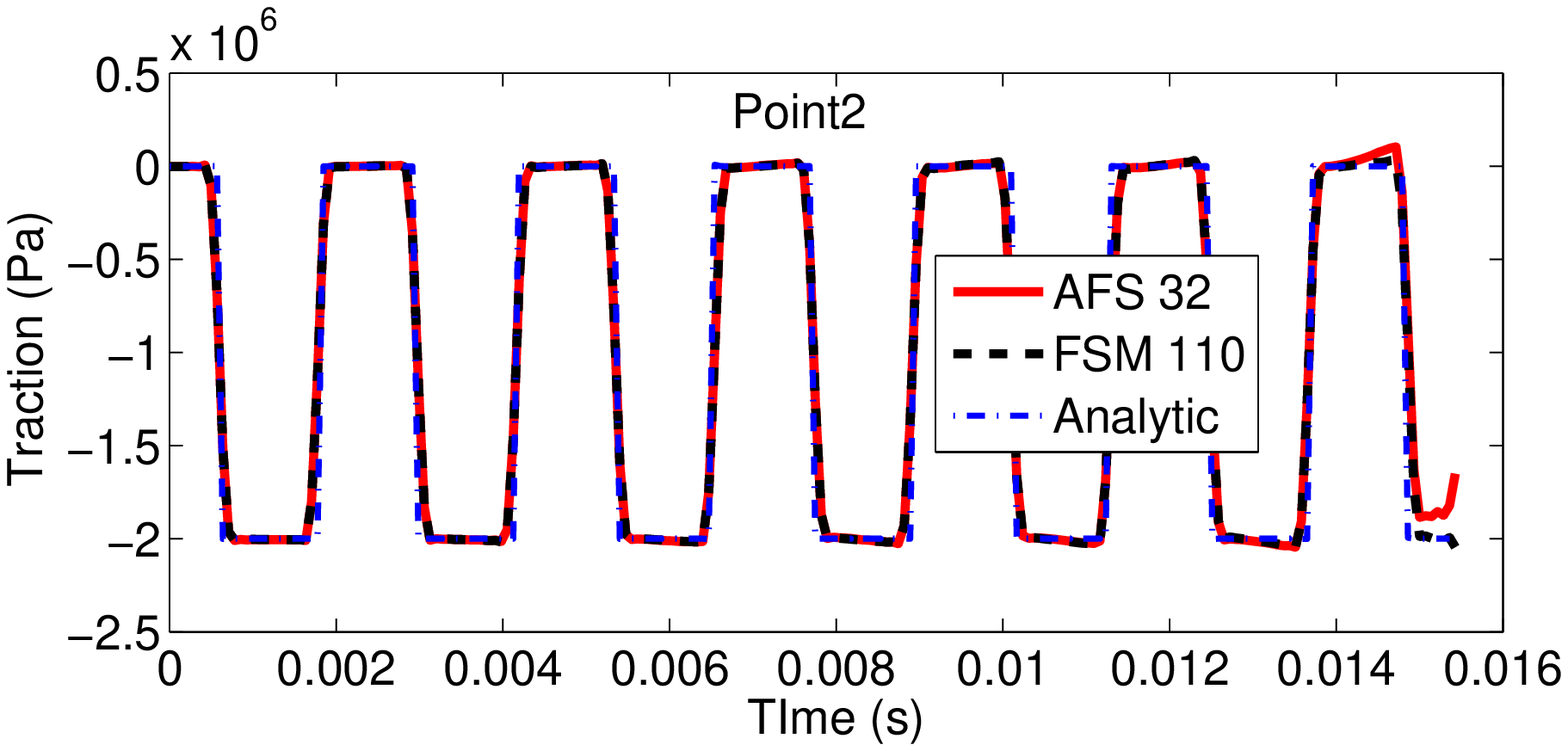,width=0.7\textwidth}
\caption{Time-domain responses of observation points $1$ and $2$.}
\label{fig-n1-td-obs12}
\end{figure}

\begin{figure}[hbt]
\centering
\epsfig{figure=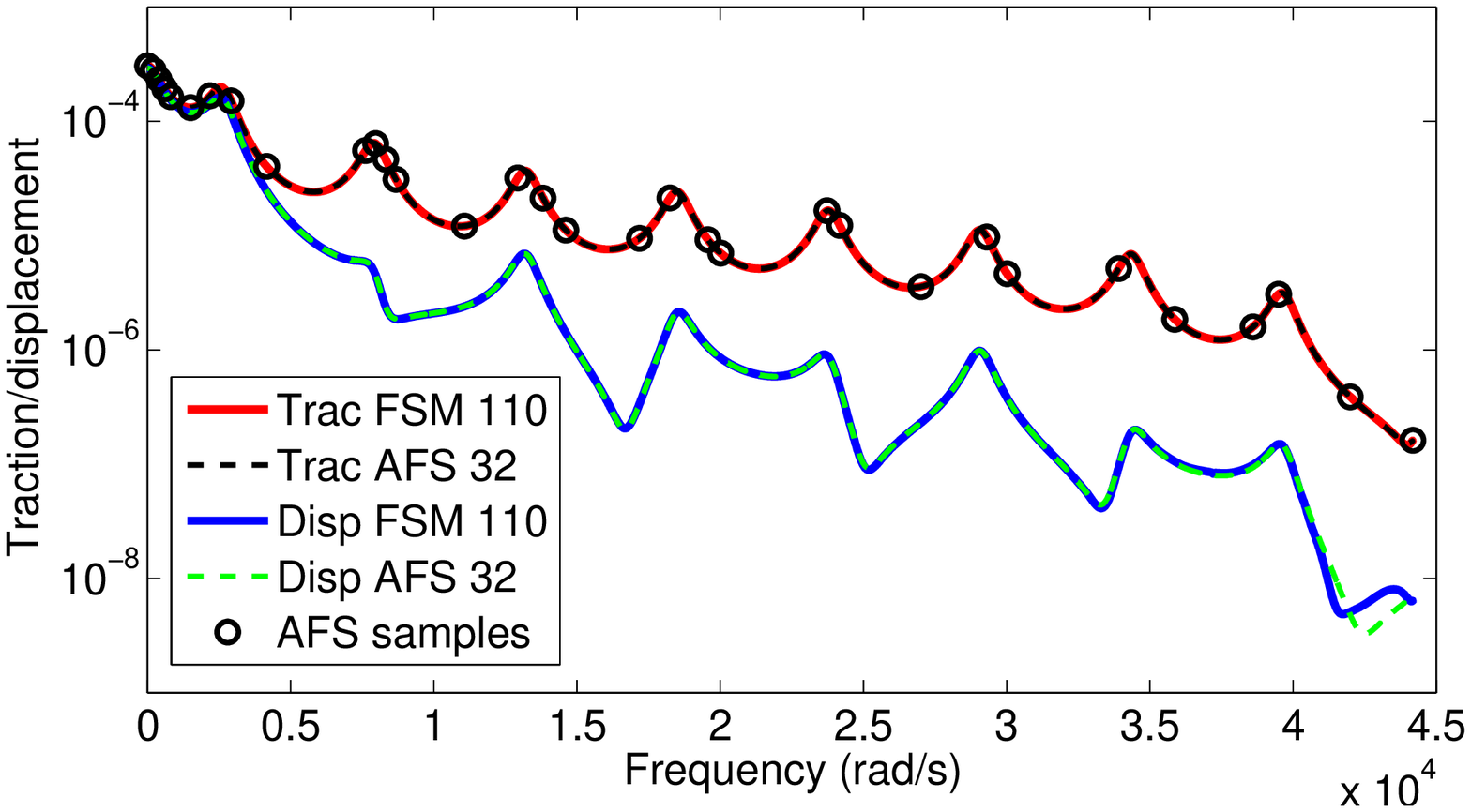,width=0.7\textwidth}
\caption{The amplitudes of the frequency-domain responses of two test points.
Note that the traction responses have been properly amplified to fit the figure. }
\label{fig-n1-fd-tst67}
\end{figure}

\begin{figure}[hbt]
\centering
\epsfig{figure=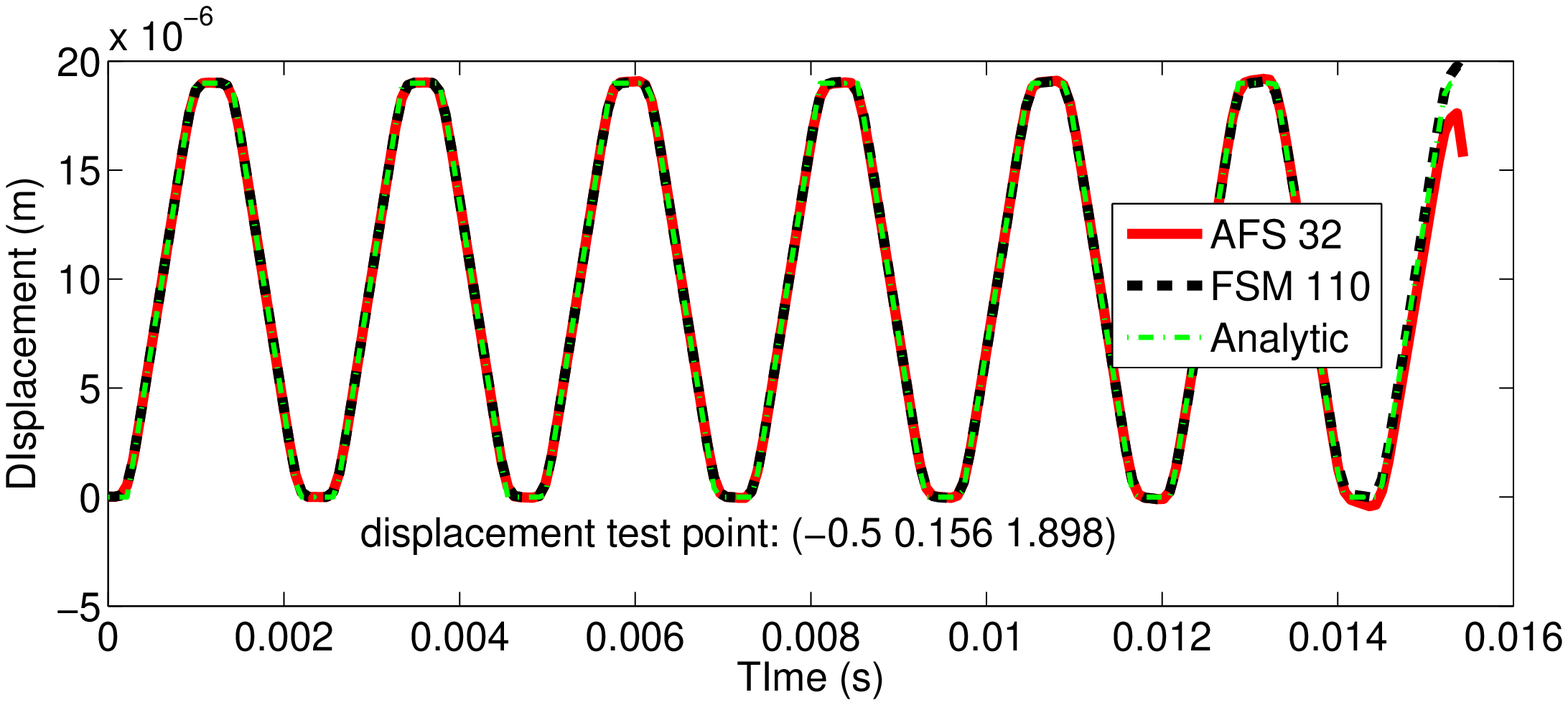,width=0.7\textwidth}
\epsfig{figure=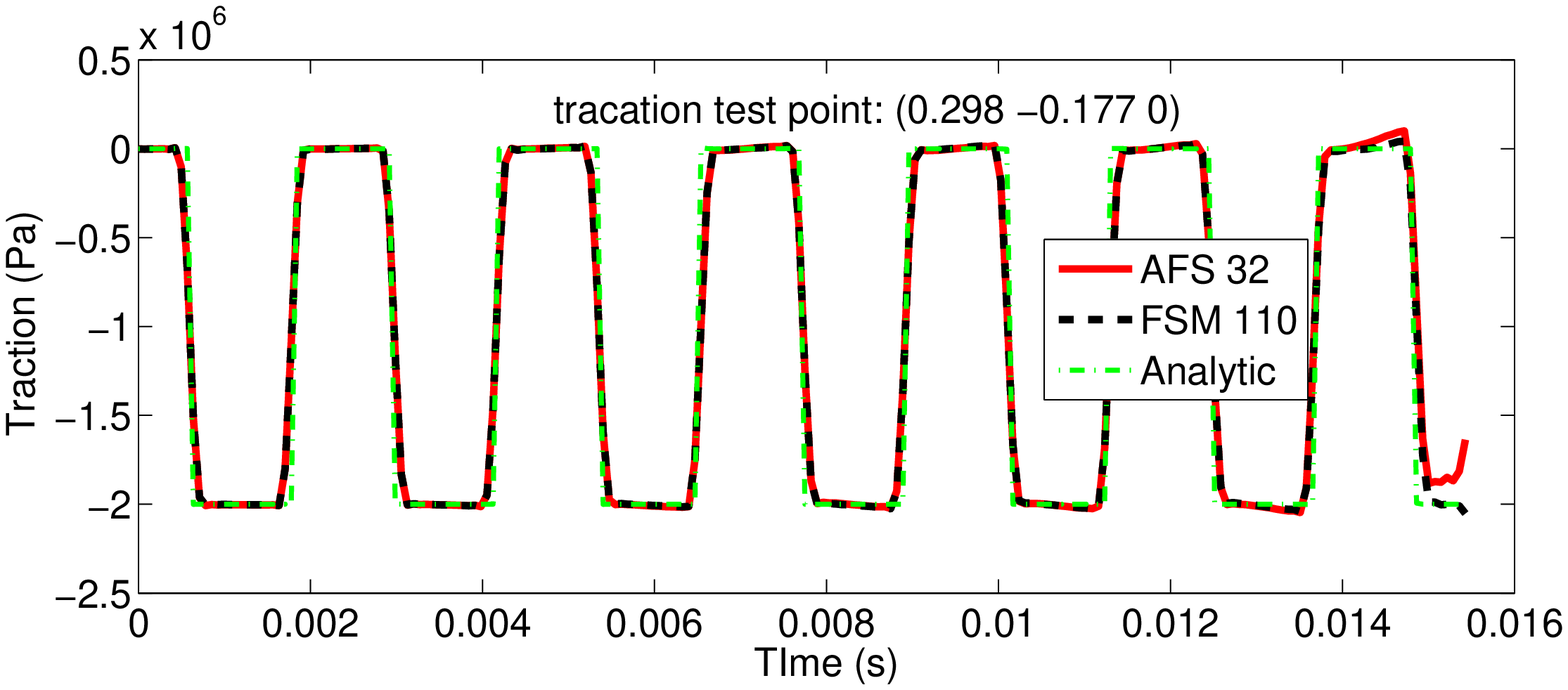,width=0.7\textwidth}
\caption{Time domain-responses of the two test points.}
\label{fig-n1-td-tst67}
\end{figure}

The AFS converges with $N_{\rm{c}}=32$. The results are compared against the FSM with $N_{\rm{c}}=110$. Figures \ref{fig-n1-fd-obs12} and \ref{fig-n1-td-obs12} illustrate the computed frequency- and time-domain responses of two observation points indexed as 1 and 2 in Table \ref{tab-NE1-vferr}. Figures \ref{fig-n1-fd-tst67} and \ref{fig-n1-td-tst67} are the similar results of two test points 6 and 7. In general, the results of AFS and FSM, both in frequency- and time-domain, coincide very well with each other and also with the analytical solutions. Here we note that the discrepancies  in Figures \ref{fig-n1-fd-obs12} and \ref{fig-n1-fd-tst67} at large frequencies are actually very small; the amplitude of the discrepancy is of 4 orders smaller than the maximal values of the corresponding spectrum. Therefore, they would not affect too much the time-domain results. For obtaining time-domain responses, and contribution of low frequencies with large amplitude is more dominant, and thus must be more accurately accounted. This is reflected in the design of error indicator \eqref{eq:ferr} in selecting new sampling frequencies, where the denominator is the maximum value of the fitted function.

\subsection{Example 2: Plate with a hole subjected to an impact loading} \label{S-S-NE2}

\begin{figure}[hbt]
\centering
\epsfig{figure=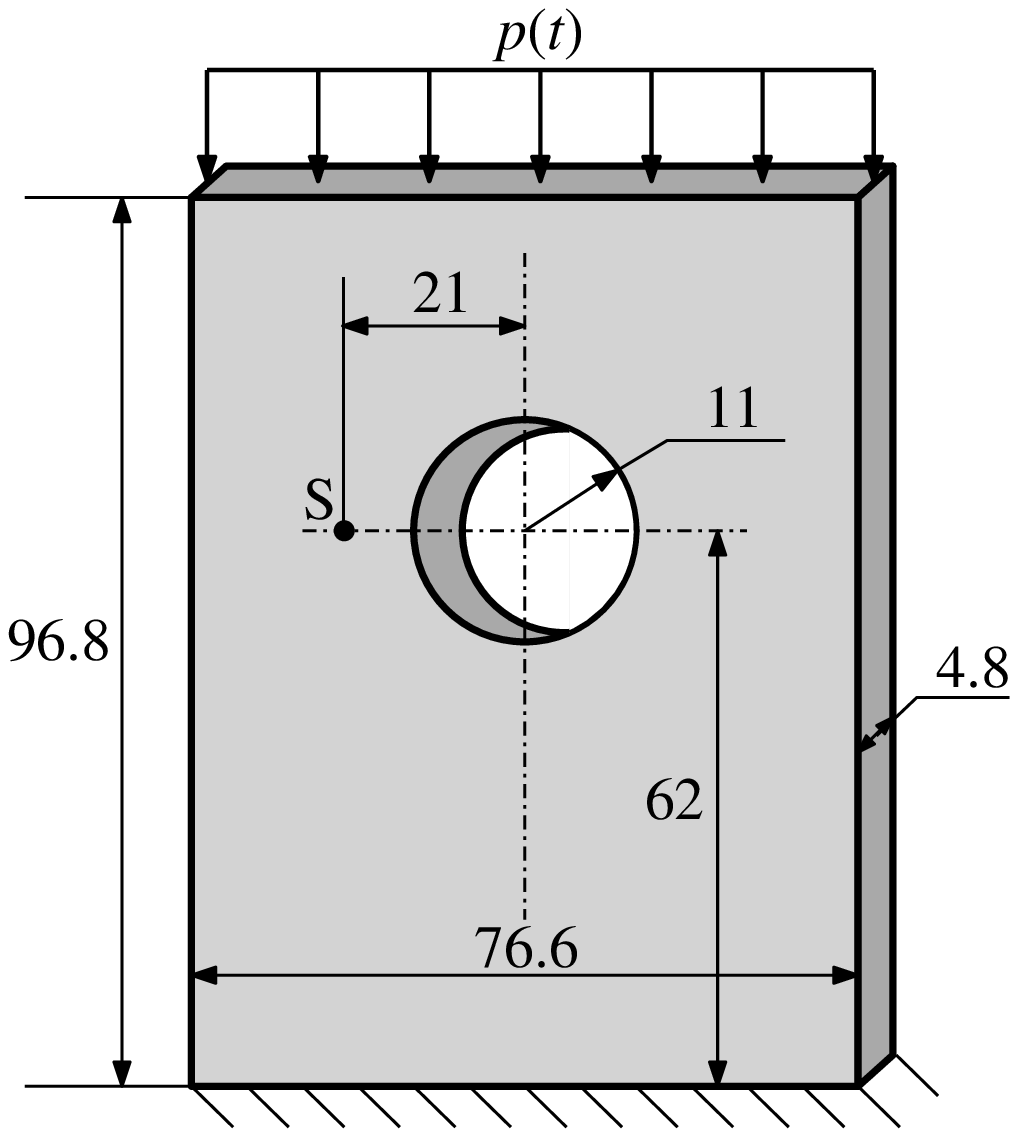,width=0.25\textwidth}\epsfig{figure=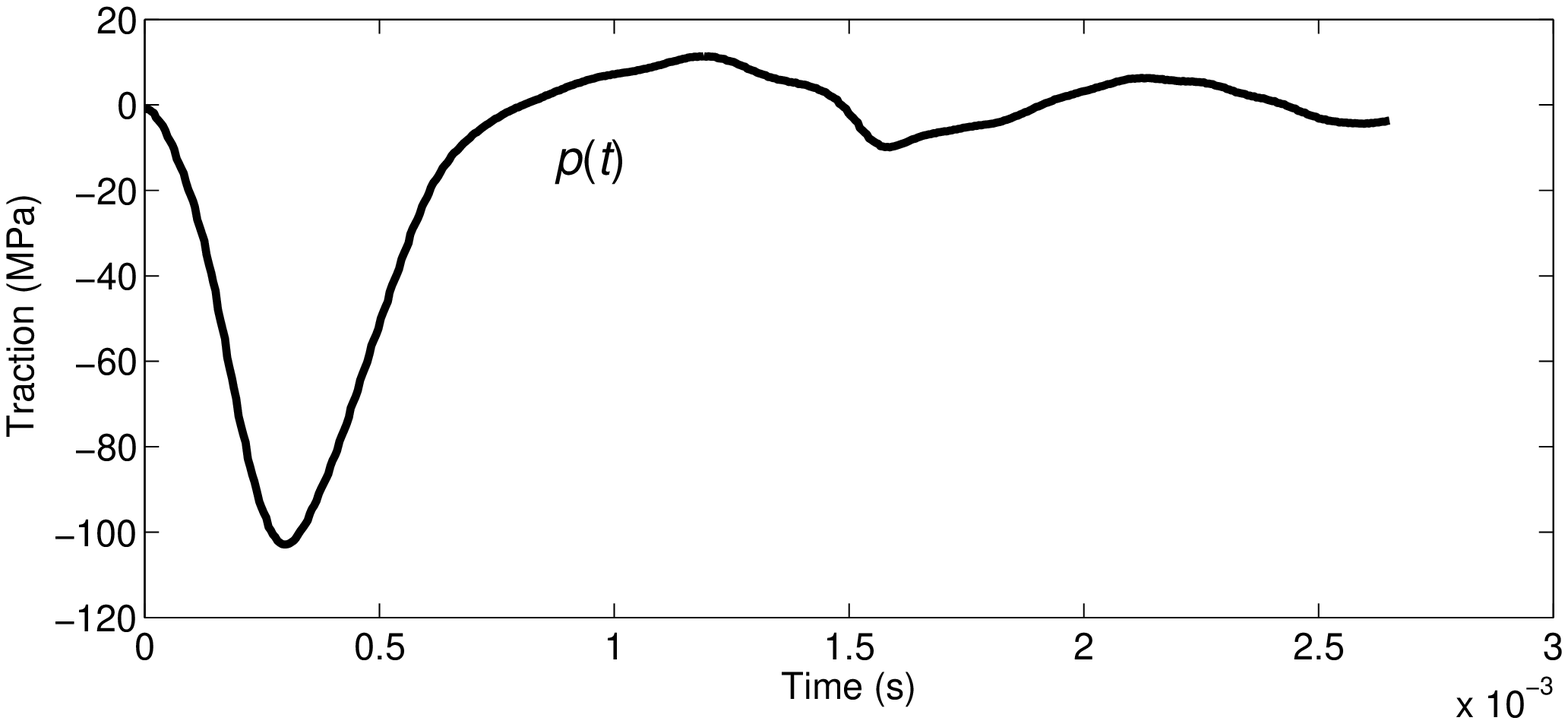,width=0.60\textwidth}
\caption{The model of the plate with central hole (unit: mm) and the loading history on the upper side.}
\label{fig-plate-hole}
\end{figure}

Here we test the proposed AFS algorithm by using a more realistic problem. The model is an aluminium plate with a central hole. The plate is fixed at the lower end and is impacted by a dynamic loading $p(t)$ on the upper end. The loading is assumed to be uniformly distributed on the upper side. The model and the load history are diagrammed in Figure
\ref{fig-plate-hole}. The material
properties of the aluminum plate are: $E=69$ GPa, $\rho=2.7 \cdot
10^3$ kg/m$^3$ and $\nu=0.3$.
The strain response at point S beside the hole is measured by a strain gauge, and is used to verify the correctness of the present results.

A boundary mesh with 2536 triangular quadratic elements is used in the frequency-domain BEM simulations. The vertical impact loading is applied to the upper side, and the lower side is fixed. All other surfaces are free of traction. The solution period is $T=0.00265$s. In FSM we use $N_{\rm s } = 160$, thus $N_{\rm c } = 81$ frequencies are solved.

\begin{figure}[hbt]
\centering
\epsfig{figure=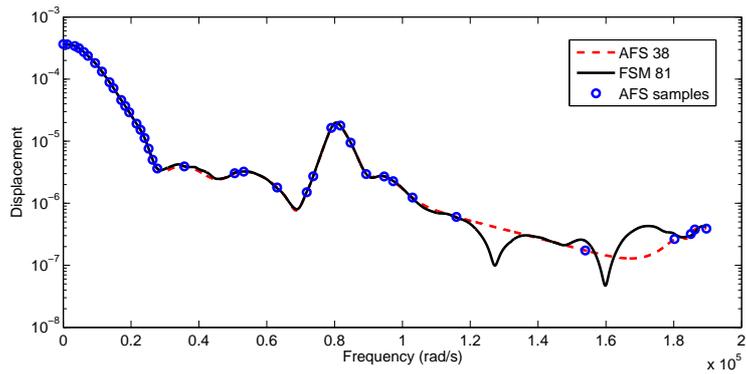,width=0.7\textwidth}
\caption{Amplitudes of the vertical displacement at Point-S in the frequency domain.}
\label{fig-plate-dispfd}
\end{figure}

\begin{figure}[hbt]
\centering
\epsfig{figure=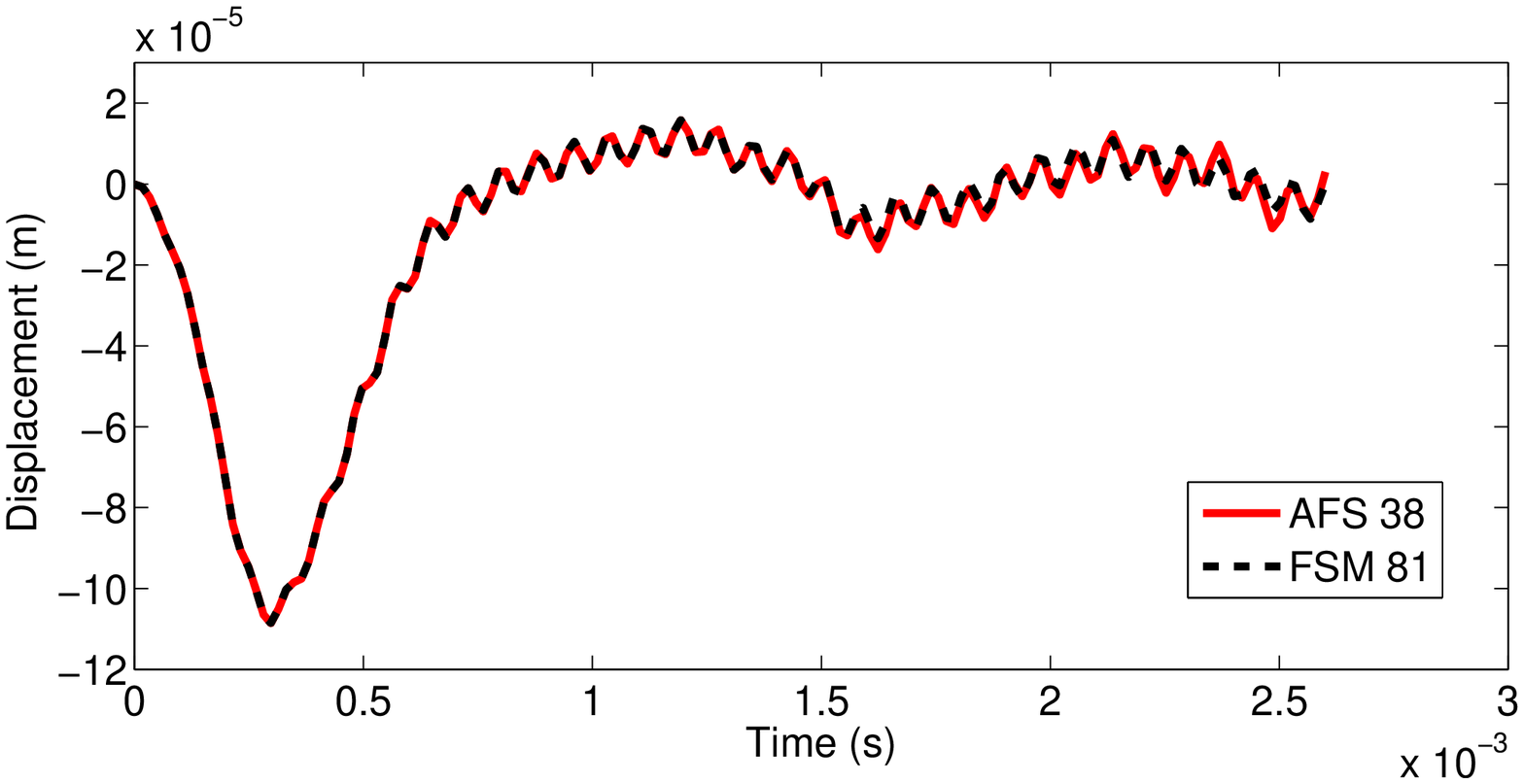,width=0.7\textwidth}
\epsfig{figure=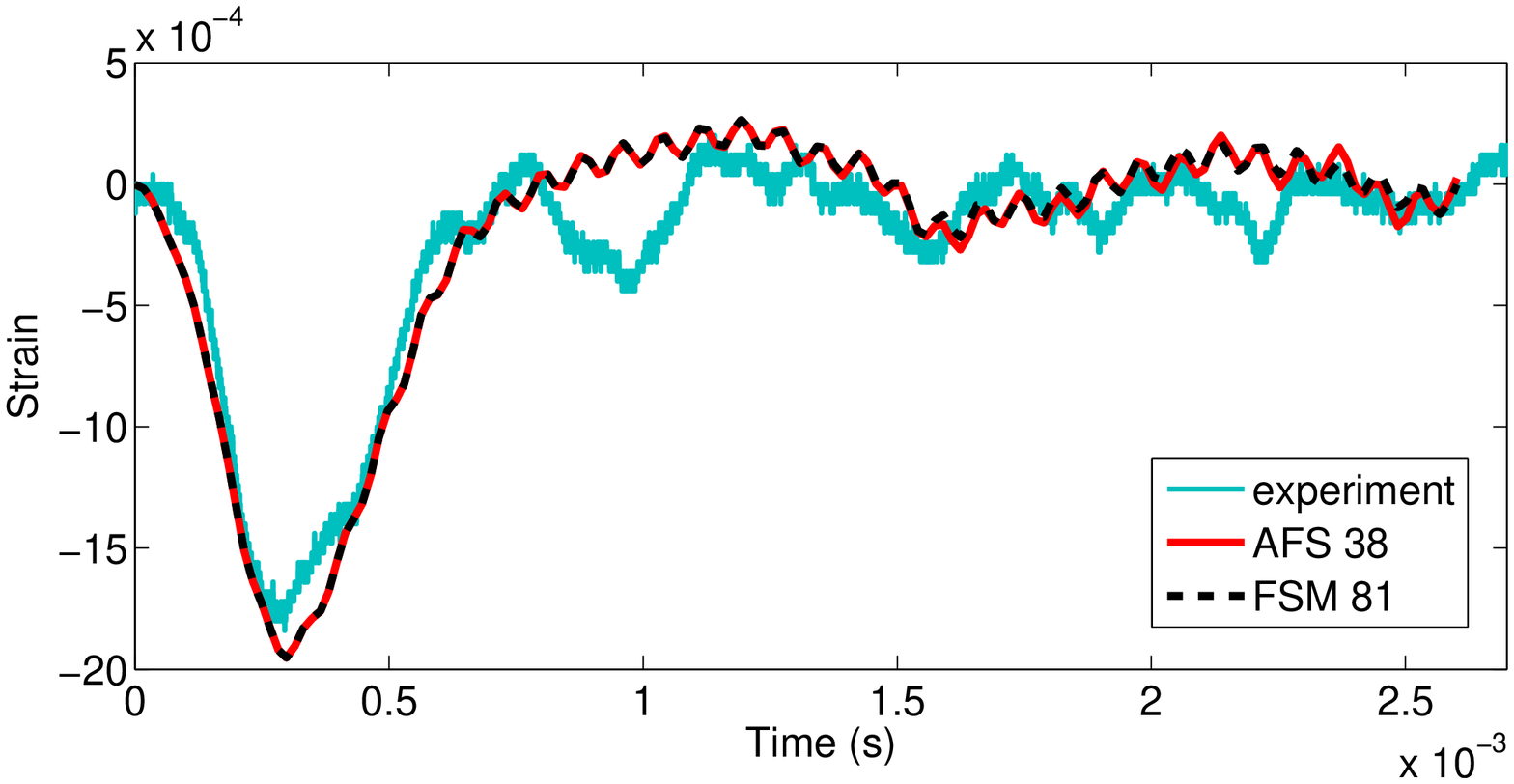,width=0.7\textwidth}
\caption{Time-domain responses of vertical displacement and strain at Point-S.}
\label{fig-plate-ds}
\end{figure}

In AFS computation, we use 5 observation responses and 2 testing responses. The other controlling parameters of the AFS are set according to the description in Section \ref{S-S-imp}. The AFS algorithm converges at $N_{\rm c } = 38$. Figure \ref{fig-plate-dispfd} illustrates the amplitudes of the vertical displacements at Point-S obtained by AFS and FSM in the frequency domain. The circles indicate the computed response of the AFS. The sampling frequencies of the AFS are more concentrated in the frequency band with large amplitude, since the accuracy of fitting in this area has the most essential influence on the time-domain responses.
The results of the two methods agree very well at frequencies with large amplitude. Discrepancy appears at frequencies with relatively small amplitude; this, however, would not affect the time-domain displacement and strain responses, as shown in Figure \ref{fig-plate-ds}, where good agreement between the two methods in the entire period can be observed.

\subsection{Example 3: Elastic block with 98 spherical cavities} \label{S-S-NE3}

\begin{figure}[hbt]
\centering
\epsfig{figure=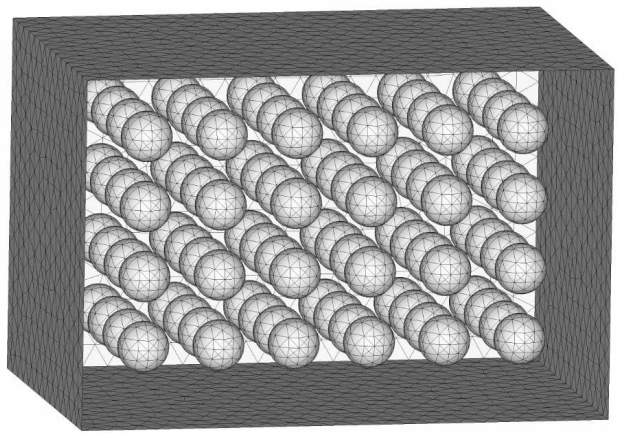,width=0.5\textwidth}
\caption{Solid block with 98 spherical cavities.}
\label{fig-prism}
\end{figure}

As the last example, the AFS is used to solve large-scale wave propagation
problem in porous solids. A model consisting of a solid block with 98 identical spherical cavities
(see Fig. \ref{fig-prism}; the front surface is removed to show the cavities) was built and simulated. The size of the block is
$1\text{m}\times 1 \text{m} \times 1.5 \text{m}$.
The spherical cavities are equally distributed inside the block
and their radius is 0.072 m.
The block is fixed at the right surface ($z=1.5$) and is subjected to a uniform step traction
$p(t)$ in $z$-direction at the left surface ($z=0$).  All other surfaces are traction free.
The material is aluminum with properties similar to those defined in Section 6.2.

In the simulation, the surfaces of the block and the cavities are
discretized into $18080$ triangular quadratic elements;
hence, the total number of unknowns is more than 0.3 million.
The response period is set to be $T=0.035$s.
In FSM, the number of sampling frequencies is set to be $478$, thus $N_{\rm{c}}=240$ frequency-domain BEM
solutions are performed.
In AFS, $6$ observation points are selected in total, with two points located at the left surface
monitoring tractions, and the other four points monitoring displacements.

\begin{figure}[hbt]
\centering
\epsfig{figure=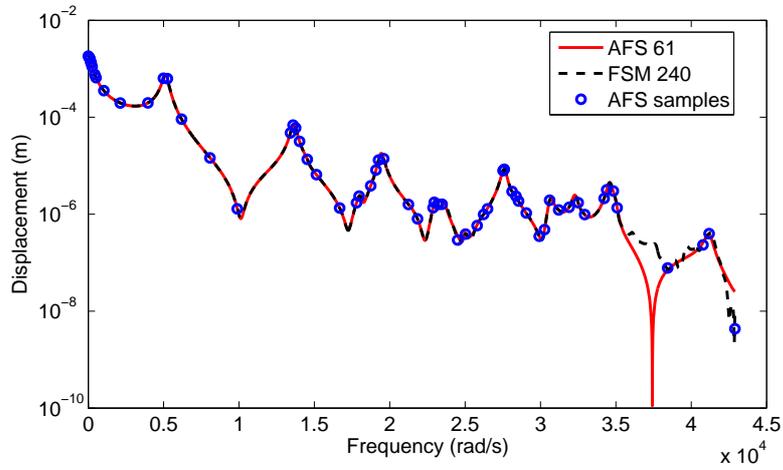,width=0.7\textwidth}
\caption{Amplitude of $z$-displacement at point A in the frequency domain.}
\label{fig-ex3dispfd}
\end{figure}

\begin{figure}[hbt]
\centering
\epsfig{figure=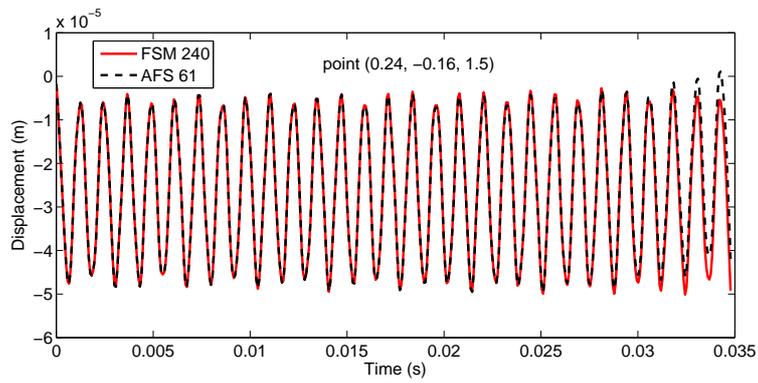,width=0.7\textwidth}
\caption{$z$-displacement at point A in the time domain.}
\label{fig-ex3disp}
\end{figure}

\begin{figure}[hbt]
\centering
\epsfig{figure=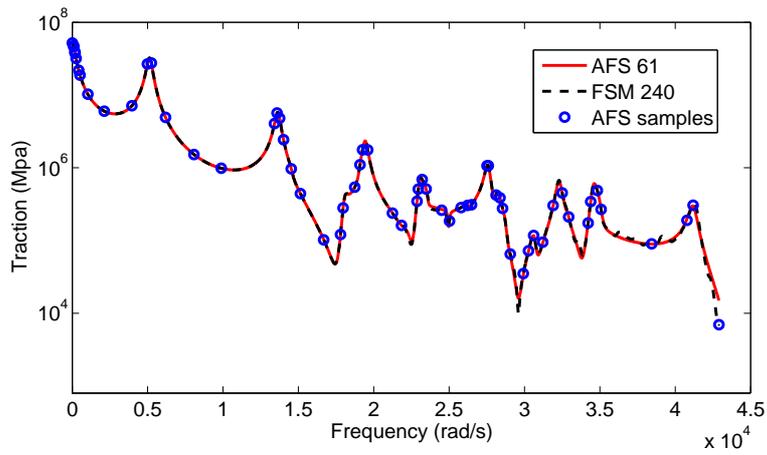,width=0.7\textwidth}
\caption{Amplitude of $z$-traction at point B in the frequency domain.}
\label{fig-ex3tracfd}
\end{figure}

\begin{figure}[hbt]
\centering
\epsfig{figure=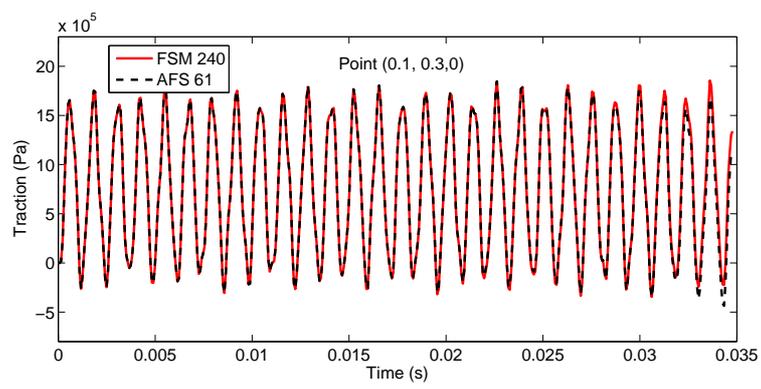,width=0.7\textwidth}
\caption{$z$-traction at point B in the time domain.}
\label{fig-ex3trac}
\end{figure}

The AFS converges at $N_c=61$. The total time consumption is about $14.9$h.
As a comparison, the entire computation of the FSM takes about $57.6$h,
which is about four times longer than the AFS.
To illustrate the accuarcy of the AFS, two non-observation points, A($0.24,\, -0.16,\, 1.5$) and B($0.1,\, 0.3,\, 0$), are selected and their responses
computed by the two methods are compared. Figures \ref{fig-ex3dispfd} and \ref{fig-ex3disp} exhibit the frequency- and time-domain displacements of Point-A in $z$-direction, respectively. Figures \ref{fig-ex3tracfd} and \ref{fig-ex3trac} exhibit the $z$-traction of Point-B in the similar way. The circles indicate the location of sampling frequencies of the AFS. They are more concentrated to the region where the frequency response is irregular and where small relative error \eqref{eq:ferr} is required. The discrepancy in Figures \ref{fig-ex3dispfd} and \ref{fig-ex3tracfd} at the end of the frequency spectral is actually 4 orders smaller than the maximal value of the respective spectrum. The time-domain solutions of the AFS and FSM are well coincided.

\section{Conclusions and discussions}\label{S-conclution}

This paper is devoted to the effective reduction of $N_{\rm{c}}$, the number of frequency-domain BEM solutions, in the frequency-domain approach for large-scale time-domain BEM simulation. We come up with an adaptive frequency sampling algorithm, i.e., the AFS, that sequentially and automatically selects the ``optimal'' sampling frequencies in the frequency-sweeping process. The AFS algorithm is inspired by two facts. First, the  frequency spectrum of real structures is featured by a series of peaks and fluctuations and it is often very smooth between two peaks. Second, low frequency spectrum contributes much more than high frequency spectrum, and thus the accuracy of the former is more important to guarantee the quality of the time-domain response. The traditional frequency-domain approach, or more specifically, the Fourier series method (FSM), uses equally-spaced sampling frequencies which might be good in computing the inverse Laplace transform, but not good in approximating the frequency spectrum of real structures. For the latter task, it is more reasonable to distribute more sampling points around the peaks and fewer ones in the smooth part. Along this way it is possible to considerably reduce $N_{\rm{c}}$, which finally leads to substantial reduction of the computational cost of the frequency-domain approach. This is verified and realized by the AFS algorithm.

Our AFS algorithm relies on the stable and accurate rational approximation of the frequency-domain responses, which is accomplished by using the VF algorithm \cite{GS99}. In the VF algorithm, the poles of the dynamic system is first extracted from the frequency-domain data at hand. These poles are then used to fit the frequency-domain responses of each DOF and to evaluate the corresponding time-domain response. In this sense, our AFS bears similarity to the well-known method of modal superposition in finite element analysis. The main difference is that here the poles obtained by the VF consist of not only some natural frequencies of the system, but also some ``fictitious frequencies'' which are necessary to guarantee the accuracy of the frequency-domain fitting in the given frequency band. We suggest to randomly select $5\sim 10$ DOFs in the pole-identification process in order to keep the computational cost of the VF algorithm affordable. We found that by using these poles one can always achieve 4-digital accuracy in fitting the frequency-domain responses of all the DOFs, and that using more DOFs in the pole-identification process does not improve too much the accuracy of the fitting.

The computational efficiency and accuracy of the AFS is illustrated and verified by the comparison with the FSM using three typical examples in Section \ref{S-ne}. We conclude that the AFS algorithm can considerably reduce $N_{\rm{c}}$ against the FSM while still maintain the accuracy of the time-domain results. A 4-times reduction of $N_{\rm{c}}$, and thus of the total computational time, is achieved in the large-scale simulation of Section \ref{S-S-NE3}.
The AFS algorithm presented in this paper can be readily adopted to accelerate the frequency-sweeping problems in structural dynamics and acoustics.

\section*{Acknowledgements}

JX gratefully acknowledges the financial supports from the National Science
Foundations of China under Grants 11102154 and 11472217, Fundamental Research Funds for the Central Universities and the Alexander von Humboldt Foundation (AvH) to support his fellowship research at the Chair of Structural Mechanics, University of Siegen, Germany.

\bibliography{fdas_bib}

\newpage

\begin{figure}[hbt]
\centering
\epsfig{figure=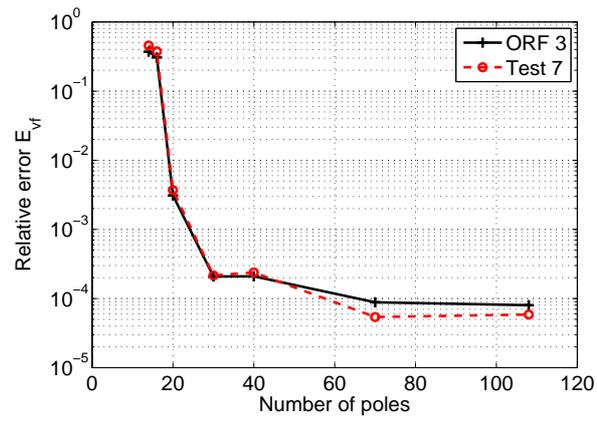,width=0.5\textwidth}
\caption{Frequency domian responses of the two test points.
Note that the traction responses have been properly amplified to fit the figure. }
\label{fig-fdtest}
\end{figure}

\end{document}